# A Nonlinear Pairwise Swapping Dynamics to Model the Selfish Rerouting Evolutionary Game

## Wen-yi Zhang, Wei Guan, Ji-hui Ma, Jun-fang Tian

**Abstract** In this paper, a nonlinear revision protocol is proposed and embedded into the traffic evolution equation of the classical proportional-switch adjustment process (PAP), developing the present nonlinear pairwise swapping dynamics (NPSD) to describe the selfish rerouting evolutionary game. It is demonstrated that i) NPSD and PAP require the same amount of network information acquisition in the route-swaps, ii) NPSD is able to prevent the over-swapping deficiency under a plausible behavior description; iii) NPSD can maintain the solution invariance, which makes the trial and error process to identify a feasible step-length in a NPSD-based swapping algorithm is unnecessary, and iv) NPSD is a rational behavior swapping process and the continuous-time NPSD is globally convergent. Using the day-to-day NPSD, a numerical example is conducted to explore the effects of the reaction sensitivity on traffic evolution and characterize the convergence of discrete-time NPSD.

## 1 Introduction

Wardrop (1952) formally put forward two route choice principles which later were termed as user equilibrium (UE) principle and system optimization (SO) principle. Later, Beckmann et al. (1956) created a convex mathematical programming model equivalent to the UE condition, stimulating rich advances in the traffic assignment theory (Sheffi 1985; Patriksson 1994; Peeta and Ziliaskopoulos 2001; Jin 2014). A central task in the traditional traffic assignment

W. -Y. Zhang (Corresponding author)
MOE Key Laboratory of Urban Transportation Complex System Theory and Technology, Beijing Jiaotong University, Beijing 100044, China
e-mail: wyzhang@bjtu.edu.cn
W. Guan
School of Traffic and Transpotation, Beijing Jiaotong University, Beijing 100044, China
e-mail: weig@bjtu.edu.cn
J. -H. Ma
School of Traffic and Transpotation, Beijing Jiaotong University, Beijing 100044, China
e-mail: jhma@bjtu.edu.cn
J. -F. Tian
College of Management and Economics, Tianjin University, Tianjin 300072, China
e-mail: jftian@tju.edu.cn

studies is to identify the static (or dynamic) equilibrium states. For the long term traffic planning, finding the equilibrium states is important and valuable. However, these studies face a reality: the equilibrium state may not be reachable; even if reachable, sufficient time and route swaps are needed. Therefore, besides the equilibrium states, traffic assignment should still attach importance to model the route-swapping dynamical process (also referred to as day-to-day traffic assignment). The day-to-day traffic modeling approach pays more attention to the disequilibrium evolution process rather than the final equilibrium state. As mentioned in Watling and Hazelton (2003), it has great flexibility which allows a wide range of behavior rules, levels of aggregation, and traffic modes to be synthesized into a uniform framework. Therefore, the day-to-day traffic modeling approach is viewed as the most appropriate one to analyze the traffic equilibrium process (He et al. 2010).

The day-to-day route-swapping models can be divided into several classes. The deterministic models (e.g., the models summarized in Yang and Zhang 2009) and stochastic models (Cantarella and Cascetta 1995; Watling 1999; Bie and Lo 2010; Xiao and Lo 2014) are classified by the network stochasticity; the path-based models (Yang and Zhang 2009; Cantarella and Cascetta 1995; Watling 1999) and link-based models (He et al. 2010; Han and Du 2012; He and Liu 2012; Guo et al. 2013; Di et al. 2014) are classified by the carrier of traffic evolution. He et al. (2010) pointed out two deficiencies in the path-based models, i.e., the path-flow-nonuniqueness problem and the path-overlapping problem. The first deficiency is a technical deficiency, which can be resolved by estimating the most likely path flow pattern in theory (Bar-Gera 2006) or tracing the traffic flows in practice. The second one is a behavioral deficiency. He et al. (2010) argued that this problem was likely to exist generally in the deterministic path-based models. Smith and Mounce (2011) proposed a splitting rate rerouting model which, however, did not exhibit the anomaly identified by He et al. (2010). The splitting rate model is the splitting rate version of the proportional-switch adjustment process (PAP) introduced by Smith (1984). Thus, this paper does not intend to address the path-overlapping problem but focus on developing a path-based route-swapping process with a sounder behavior realism under the pairwise swapping framework. We argue that, from the behavioral perspective, a path-based rerouting model is more fundamental in modeling the individual route-swapping behaviors since it is the route-swapping rather than link-swapping which drives the traffic evolution. According to the focus of this study, here we just present a brief review on the deterministic path-based rerouting models.

Except for the aforementioned PAP (Smith 1984), there are many other deterministic path-based rerouting models, for instance, the simplex gravity flow dynamics (Smith 1983), the network tatonnement process (Friesz et al. 1994), the projected dynamic system (Zhang and Nagurney 1996; Nagurney and Zhang 1997), and the BNN swapping process (Yang 2005). Yang and Zhang (2009) proved the five rerouting dynamics are rational behavior adjustment processes (RBAP) whose stationary link flows are UE. Among the abovementioned five RBAPs, PAP is the most natural and has the simplest formulation. PAP has stimulated various extensional applications (e.g., Smith and Wisten 1995; Huang and Lam 2002; Peeta and Yang 2003; Yang et al. 2007; Mounce and Carey 2011; Smith and Mounce 2011, etc.). Based on PAP, Cho and Hwang (2005a, 2005b) proposed a stimulus-reaction dynamic model. Li et al. (2012) developed an excess-travel-cost PAP. From a new perspective, Jin (2007) proposed a J-dynamic model whose stationary state was not equivalent to UE. Therefore, the J-dynamic model is not a RBAP. Considering that RBAP has been viewed as a basic behavior rule to model the rerouting behavior, this study only discusses the RBAP-based models.

PAP is initially developed to analyze the stability of UE. The essence of PAP is the traffic evolution equation, where a linear revision protocol is embedded to model the individual rerouting mechanism. The traffic evolution equation of PAP can provide a good behavior approximation for the Markov evolution games, and the majority of the existing population evolution dynamics can be rewritten in the PAP form with a suitable revision protocol (Hofbauer 2011 and Sandholm 2011). However, PAP's revision protocol assumes the travelers adjust their routes day after day in a linear pairwise way; this does not quite accord with the fact that human's behaviors are more likely to be nonlinear. In order to capture the rerouting behaviors better, this paper proposes a nonlinear revision protocol and embeds it into the traffic evolution equation of PAP, deriving the present nonlinear pairwise swapping dynamics (NPSD). NPSD can not only provide a plausible description for travelers' rerouting behaviors but also add a new algorithmic device to solve the traffic equilibrium.

The remaining context is organized as follows. In Section 2, PAP is given and analyzed. Section 3 elaborates on the proposed NPSD model. Section 4 discusses the stability of NPSD. Applying the day-to-day NPSD, a numerical example is performed in Section 5 to explore the impacts of reaction sensitivity on the network traffic evolution and stability of discrete-time NPSD (DNPSD). Section 6 concludes the whole study and suggests some future works.

## 2 PAP

The continuous-time PAP (CPAP), firstly proposed by Smith (1984), is formulated as follows:

$$\dot{f}_k^{wt} = \sum_{p \in K^w} f_p^{wt} \rho_{pk}^{wt} - f_k^{wt} \sum_{p \in K^w} \rho_{kp}^{wt} \quad \forall k \in K^w, w \in W, t \leq T, \tag{1}$$

where

$$\rho_{kp}^{wt} = \kappa \max\left(0, C_k^{wt} - C_p^{wt}\right). \tag{2}$$

Here, $T$ is the terminal time of swapping, $t$ is the time index, $W$ is the set of OD-pairs, $w$ is the OD-pair index, $p$ and $k$ are the route indices, $K^w$ is the set of effective routes between OD-pair $w$, $f_k^{wt}$ is the flow on route $k$ between OD-pair $w$ at time $t$, $\dot{f}_k^{wt}$ is the derivative of $f_k^{wt}$ with respect to $t$, $C_k^{wt}$ is the cost on route $k$ between $w$ at time $t$, $\kappa$ is a small positive constant, and $\rho_{kp}^{wt}$ is the revision protocol to estimate the proportion of flow that swaps from $k$ to $p$ between $w$ at time $t$.

Eq. (1) is the traffic evolution equation which reflects a fact that the change of traffic flow on a route equals to the total flow swapping onto that route minus the total flow swapping off it. On the right side of Eq. (1), the first expression formulates the total flow swapping onto route $k$ and the second one formulates the total flow swapping off it. Replacing $\dot{f}_k^{wt}$ by $f_k^{w(t+1)} - f_k^{wt}$ in Eq. (1), the discrete-time PAP (DPAP) is derived.

According to Eq. (2), we can conclude $\rho_{kp}^{wt} \geq 0$ and further

$$\rho_{kp}^{wt} \begin{cases} = 0 & \text{if } C_k^{wt} \leq C_p^{wt} \\ = \kappa\left(C_k^{wt} - C_p^{wt}\right) > 0 & \text{if } C_k^{wt} > C_p^{wt} \end{cases} \quad \forall k, p \in K^w, w, t. \tag{3}$$

Eq. (3) implies that i) traffic only swaps from a costly route to the less costly ones in the same OD-pair, and ii) the swapping ratio is proportional to the absolute cost difference between two routes. Obviously, $\rho_{kp}^{wt} \leq 1 \ \forall k, p, w, t$ is required; otherwise, an illogical consequence named over-swapping will happen and the negative route flows can be produced. The definition of over-swapping is given below.

**Definition 1** (Over-swapping). During a route-swapping, the total flow swapping off a route is larger than the initial flow on it, mathematically, it means that there exists at least a $k, w, t$ to make $\sum_p \rho_{kp}^{wt} > 1$.

Generally, a small $\kappa$ is beneficial to prevent over-swapping. Smith and Wisten (1995) presented an upper bound of such a small $\kappa$ that can avoid over-swapping, i.e., $\kappa \leq (BM)^{-1}$, where $B = \max\{C_k^{wt} | k \in K^w, w \in W, t \leq T\}$ and $M$ is the number of effective routes in the network. Obviously, this upper bound does not take much behavior realism into account; it is more like a mathematical construction. In the swapping algorithms (e.g., Nie 2003), several trial and error techniques (e.g., backtracking, Armijo search, etc.) are employed to indentify an appropriate $\kappa$ in each iteration.

He et al. (2010) gave another expression for $\kappa$ in PAP, i.e.,

$$\kappa = \frac{1}{\sum_p \sum_k \max\left(0, C_p^{wt} - C_k^{wt}\right) + H}. \tag{4}$$

Here, $H > 0$ is a reluctance parameter, and more travelers prefer maintaining the previous choices when a larger $H$ appears. Mathematically, over-swapping can be prevented by the $\kappa$ in Eq. (4). However, this formulation is still unable to provide a satisfactory behavior explanation. For instance, the denominator of Eq. (4) implies that the swapping proportion of a certain route to its candidate also depends on the pairwise cost differences excluding this route. To describe the problem more clearly, we introduce the following example network (see Fig. 1), and assume, in the beginning, that i) there is traffic on every route, and ii) Route 1 is the most costly, next are Route 2, Route 3 and Route 4, respectively.

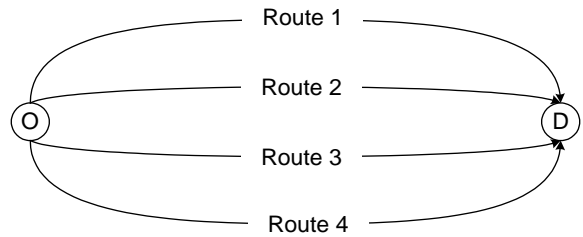

**Fig. 1** The example network

According to the PAP with $\kappa$ formulated by Eq. (4), traffic can swap from Route 1 to Route 3, and the swapping ratio also depends on the pairwise cost differences of irrelative routes, e.g., Route 2 and Route 4. Obviously, it is unreasonable.

It is known from above analysis that a good revision protocol for PAP should not only prevent the over-swapping but also possess a plausible behavior description. To this end, this paper suggests a new revision protocol and embeds it into the traffic flow evolution equation of PAP, deriving the present NPSD model.

## 3 Nonlinear pairwise route-swapping dynamics

Firstly, we present four assumptions for the present route-swapping dynamic model. These assumptions are also applied in PAP (Smith 1984).

1) A traveler only changes his route to a less costly one, and, at least some travelers, if not all, will do so unless all the travelers were all on the least cost routes the previous day.
2) The travel demand is inelastic.
3) Every driver has perfect information on the previous and current traffic network.
4) Swap decision is based on the travel experience of the previous day.

Based on the above four assumptions, the continuous-time NPSD (CNPSD) comprises an evolution equation, a nonlinear revision protocol and a feasible initial condition which are formulated by Eq. (1), Eq. (5) and Eq. (7), respectively. Replacing $\dot{f}_k^{wt}$ by $f_k^{w(t+1)} - f_k^{wt}$ in Eq. (1), the discrete-time NPSD (DNPSD) is derived.

$$\rho_{kp}^{wt} = \begin{cases} \max\left(0, \dfrac{1}{|R_k^{wt}|}\left(1-\exp\left(-\theta^w\left(C_k^{wt}-C_p^{wt}\right)\right)\right)\right), & \text{if } R_k^{wt} \neq \varnothing \\ 0, & \text{if } R_k^{wt} = \varnothing \end{cases} \quad \forall k, p, w, t \qquad (5)$$

with

$$R_k^{wt} = \left\{p \in K^w \middle| C_p^{wt} < C_k^{wt}\right\} \quad \forall k, w, t \qquad (6)$$

and the initial route flows satisfy

$$\sum_k f_k^{w0} = d^w \ \forall w; f_k^{w0} \geq 0 \ \forall k, w. \qquad (7)$$

Here, $R_k^{wt}$ is the set of candidate routes for route $k$ between OD-pair $w$ at time $t$, $|R_k^{wt}|$

is the number of routes in $R_k^{wt}$, $d^w$ is the travel demand between $w$, $\theta^w$ is a positive constant that reflects the travelers' reaction sensitivity between $w$, and the travelers react more sensitively as $\theta^w$ increases.

The revision protocol $\rho_{kp}^{wt}$ in Eq. (5) is used to estimate the proportion of traffic flow swapping from route $k$ to its candidates. Actually, we can also interpret the meaning of the revision protocol from a microscopic perspective, i.e., the possibility for a person swapping from a route to another at a time. When the amount of route flow is large enough, the two meanings can be unified according to the weak law of large numbers. Eq. (6) shows that the candidates of route $k$ should belong to the same OD-pair with route $k$ and have a strictly lower travel cost. In Eq. (7), the first equation is the equilibrium constraint of travel demand; the second one promises the path flows are nonnegative. Usually, a defined constraint for the path and link flows is also introduced, i.e., $v_a^t = \sum_w \sum_k f_k^{wt} \delta_{ak}^w \ \forall a,t$, where $v_a^t$ is the flow on link $a$ at day $t$, and 0-1 indicator $\delta_{ak}^w = 1$ if link $a$ is used by the route $k$ between $w$; otherwise, $\delta_{ak}^w = 0$.

Next, we elaborate on the behavior mechanism of the second expression on the right hand of Eq. (5) when $R_k^{wt} \neq \varnothing$. To simplify the elaboration, we interpret it from a microscopic perspective. When $R_k^{wt} \neq \varnothing$, $1-\exp\left(-\theta^w\left(C_k^{wt} - C_p^{wt}\right)\right)$ can be considered to formulate the monopolized probability of a driver swaps from route $k$ to $p$ given that $p$ is the unique candidate of route $k$. When the candidate routes are not unique, a person needs to distribute the probability on the candidates properly since one cannot choose all the candidate routes simultaneously. Here $\left|R_k^{wt}\right|^{-1}$ is designed to formulate the distributing mechanism of these monopolized probabilities; it implies that a person will reassign the monopolized probability towards each candidate by dividing the number of candidates.

Although we have assumed in prior that the network information is perfect, actually only the information of being-used route and the candidates is needed for the NPSD process, which means that NPSD does not increase the network information acquisition than PAP. Next, we state three important properties of NPSD.

**Property 1** (contrary sign). For NPSD (In the remaining context, a property or a theorem that holds for NPSD means it holds for both CNPSD and DNPSD), the following inequality holds in general, i.e.,

$$\rho_{kp}^{wt}\left(C_k^{wt} - C_p^{wt}\right) \geq 0 \ \forall k, p, w, t \tag{8}$$

**Proof.** According to Eqs. (5)-(6), it can be concluded that $\rho_{kp}^{wt} = 0$ if and only if $C_k^{wt} \leq C_p^{wt}$ and $\rho_{kp}^{wt} > 0$ if and only if $C_k^{wt} > C_p^{wt}$. Then, we have $\rho_{kp}^{wt}\left(C_k^{wt} - C_p^{wt}\right) \geq 0$. □

**Property 2** (non-over-swapping). For NPSD, during a swapping, the total flow swapping off a route cannot spill its initial value, mathematically, i.e.,

$$0 \leq \sum_p \rho_{kp}^{wt} \leq 1 \ \forall k, w, t. \tag{9}$$

**Proof.** Since $1 - \exp\left(-\theta^w \left(C_k^{wt} - C_p^{wt}\right)\right) \leq 1$, then $0 \leq \max\left(0, 1 - \exp\left(-\theta^w \left(C_k^{wt} - C_p^{wt}\right)\right)\right) \leq 1$ and $\sum_p \rho_{kp}^{wt} = \sum_{p \in K^w} \rho_{kp}^{wt} = \sum_{p \in R_k^{wt}} \left|R_k^{wt}\right|^{-1} \max\left(0, 1 - \exp\left(-\theta^w \left(C_k^{wt} - C_p^{wt}\right)\right)\right)$ further, which is less than $\sum_{p \in R_k^{wt}} \left|R_k^{wt}\right|^{-1} \times 1 = 1$ if $R_k^{wt} \neq \emptyset$ and equals to 0 otherwise. Recall $\rho_{kp}^{wt} \geq 0$, we have $\sum_p \rho_{kp}^{wt} \geq 0$, then Property 2 holds. □

**Property 3** (solution invariance). For NPSD, if the initial path flow pattern is feasible, so are the remaining path flow patterns.

**Proof.** Given a feasible initial route flow pattern, Property 3 requires that the remaining route flows are nonnegative and the travel demands are conservative, i.e.,

$$\sum_k f_k^{wt} = d^w \ \forall w; f_k^{wt} \geq 0 \ \forall k, w, t. \tag{10}$$

Eq. (10) is the feasible route flow set $\Omega$. Define a small real number $\tau > 0$, firstly, we prove that CNPSD possesses Property 3.

**Non-negativity.** According to Eq. (5) and Property 2, $\rho_{kp}^{wt} \geq 0$ and $\sum_p \rho_{kp}^{wt} \leq 1$ are guaranteed. Recall Eq. (1), we have $f_k^{w(t+\tau)} \geq f_k^{wt} + \sum_p f_p^{wt} \rho_{pk}^{wt} - f_k^{wt} = \sum_p f_p^{wt} \rho_{pk}^{wt} \geq 0$ if $f_k^{wt} \geq 0 \ \forall k, w$. In other words, $f_k^{w(t+\tau)} \geq 0$ if $f_k^{wt} \geq 0 \ \forall k, w, t$. By recurrence, it can be concluded that $f_k^{wt} \geq 0 \ \forall k, w, t$ if $f_k^{w0} \geq 0 \ \forall k, w$.

**Conservation.** Since

$$\sum_k \sum_p f_p^{wt} \rho_{pk}^{wt} = \sum_p \sum_k f_k^{wt} \rho_{kp}^{wt} = \sum_k \sum_p f_k^{wt} \rho_{kp}^{wt} \ \forall w, t,$$

according to Eq. (1), we have

$$\sum_k f_k^{w(t+\tau)} = \sum_k f_k^{wt} + \sum_k \left( \sum_p f_p^{wt} \rho_{pk}^{wt} - f_k^{wt} \sum_p \rho_{kp}^{wt} \right)$$
$$= \sum_k f_k^{wt} + 0 = \cdots = \sum_k f_k^{w0} = d^w.$$

For the above proof, let $\tau$ and $t$ take integers, the solution invariance will hold for DNPSD. Accordingly, Property 3 is proved. □

Property 3 is significant for developing a NPSD-based swapping algorithm to solve the traffic equilibrium problem. It means that the trial and error process for a feasible step-length can be omitted in a NPSD-based swapping algorithm, which helps to raise the efficiency of a swapping algorithm.

**Definition 2** (Wardrop user equilibrium (Wardrop 1952)). At the equilibrium state, all the used paths between the same OD-pair have the minimal and equal travel costs when the travel demand is fixed. Mathematically, the following expression holds at time $t$.

$$C_k^{wt} \begin{cases} = \pi^w & \text{if } f_k^{wt} > 0 \\ \geq \pi^w & \text{if } f_k^{wt} = 0 \end{cases} \forall k, w, \quad (11)$$

where $\pi^w$ is the minimum route cost between OD-pair $w$ at time $t$.

**Definition 3** (stationary path flow pattern). The stationary path flow pattern of NPSD is a set of network path flow states; starting from these path flow states, NPSD reproduces them.

**Corollary 1.** For NPSD, $\mathbf{f^t}$ (groups $f_k^{wt} \forall k, w$) belongs to the stationary path flow pattern if and only if $\dot{f}_k^{wt} = 0$ (or $f_k^{w(t+1)} = f_k^{wt}$) $\forall k, w$.

**Proof.** Since NPSD is a determined one-to-one dynamics (i.e., a determined input produces a single determined output), the necessity of Corollary 1 can be easily proved by recurrence. The sufficiency can be concluded by Definition 3. □

Subsequently, we present the equivalent relationship between stationary path flow pattern of NPSD and Wardrop user equilibrium.

**Theorem 1.** The stationary path flow pattern of NPSD is equivalent to Wardrop user equilibrium.

**Proof.** Firstly, we prove that Theorem 1 holds for CNPSD. For this, based on Corollary 1, we only need to prove that $\mathbf{f^t}$ with $\dot{f}_k^{wt} = 0 \ \forall k, w$ is equivalent to Wardrop user equilibrium.

**Sufficiency.** Suppose route $k$ has the minimum cost among the routes between $w$ (i.e., $C_k^{wt} \leq C_p^{wt} \ \forall p \in K^w$ and $f_k^{wt} > 0$), $\rho_{kp}^{wt} = 0 \ \forall p \in K^w$ is promised. According to Eq. (1), we have $\sum_p f_p^{wt} \rho_{pk}^{wt} = 0 \ \forall k, w$ if $\dot{f}_k^{wt} = 0 \ \forall k, w$. Since $\rho_{pk}^{wt}$ and $f_k^{wt}$ are both non-negative further, $f_p^{wt} \rho_{pk}^{wt} = 0 \ \forall k, w$ can be deduced. Note that $f_p^{wt} \rho_{pk}^{wt} = 0 \ \forall k, w$ implies $C_k^{wt} \leq C_p^{wt}$ $\forall f_p^{wt} \geq 0$ and $C_k^{wt} = C_p^{wt}$ if $f_p^{wt} > 0$. Denote $\pi^w = C_k^{wt}$, it means that $C_p^{wt} = \pi^w$ if $f_p^{wt} > 0$ and $C_p^{wt} \geq \pi^w$ if $f_p^{wt} = 0$, which is exactly the Wardrop user equilibrium.

**Necessity.** Suppose $\mathbf{f^t}$ is a Wardrop user equilibrium, i.e., $C_k^{wt} = \pi^w$ if $f_k^{wt} > 0$ and $C_k^{wt} \geq \pi^w$ if $f_k^{wt} = 0$, we have $C_k^{wt} \leq C_p^{wt} \ \forall p \in K^w$ if $f_k^{wt} > 0$ and $C_k^{wt} = C_p^{wt}$ if $f_k^{wt} > 0$ and $f_p^{wt} > 0$. From Eq. (5), $\rho_{kp}^{wt} = 0 \ \forall p \in K^w$ if $f_k^{wt} > 0$ and $\rho_{pk}^{wt} = 0 \ \forall p \in K^w$ if $f_k^{wt} > 0$ and $f_p^{wt} > 0$ are promised. Then we have $f_k^{wt} \sum_p \rho_{kp}^{wt} = 0 \ \forall k, w$ and $\sum_p f_p^{wt} \rho_{pk}^{wt} = 0 \ \forall k, w$ if $f_k^{wt} > 0$, which leads to $\dot{f}_k^{wt} = 0 \ \forall k, w$ if $f_k^{wt} > 0$ from Eq. (1). When $f_k^{wt} = 0$, we have $f_k^{wt} \sum_p \rho_{kp}^{wt} = 0 \ \forall k, w$, and $\rho_{pk}^{wt} = 0 \ \forall p \in K^w$ if $f_p^{wt} > 0$, then $\sum_p f_p^{wt} \rho_{pk}^{wt} = 0 \ \forall k, w$ can be also derived. Recall Eq. (1) again, we can conclude $\dot{f}_k^{wt} = 0 \ \forall k, w$ if $f_k^{wt} = 0$. Therefore, $\mathbf{f^t}$ with $\dot{f}_k^{wt} = 0 \ \forall k, w$ is at Wardrop user equilibrium.

Accordingly, Theorem 1 holds for CNPSD. Replacing $\dot{f}_k^{wt}$ by $f_k^{w(t+1)} - f_k^{wt}$ above, we obtain the proof for DNPSD. □

**Definition 4** (rational behavior adjustment process, RBAP (Yang and Zhang 2009)). A day-to-day route choice adjustment process is called a RBAP with fixed travel demand if the aggregated travel cost of the entire network decreases based on the previous day's path travel costs when path flows change from day to day. Moreover, if the path flows become stationary over days, then it is equivalent to the user equilibrium path flow.

**Theorem 2.** NPSD is a RBAP.

**Proof.** Firstly, we prove CNPSD is a RBAP. Based on Theorem 1, to prove Theorem 2, we only need to prove $\sum_w \sum_k C_k^{wt} \dot{f}_k^{wt} \leq 0$. For this, expanding $\sum_w \sum_k C_k^{wt} \dot{f}_k^{wt}$ as follows:

$$\begin{aligned}
\sum_w \sum_k C_k^{wt} \dot{f}_k^{wt} &= \sum_w \sum_k C_k^{wt} \left( \sum_p f_p^{wt} \rho_{pk}^{wt} - \sum_p f_k^{wt} \rho_{kp}^{wt} \right) \\
&= \sum_w \sum_k C_k^{wt} \left( \sum_p f_p^{wt} \rho_{pk}^{wt} \right) - \sum_w \sum_k C_k^{wt} \left( \sum_p f_k^{wt} \rho_{kp}^{wt} \right) \\
&= \sum_w \sum_k \sum_p C_k^{wt} f_p^{wt} \rho_{pk}^{wt} - \sum_w \sum_k \sum_p C_k^{wt} f_k^{wt} \rho_{kp}^{wt} \\
&= \sum_w \sum_k \sum_p C_k^{wt} f_p^{wt} \rho_{pk}^{wt} - \sum_w \sum_p \sum_k C_p^{wt} f_p^{wt} \rho_{pk}^{wt} \\
&= \sum_w \sum_k \sum_p C_k^{wt} f_p^{wt} \rho_{pk}^{wt} - \sum_w \sum_k \sum_p C_p^{wt} f_p^{wt} \rho_{pk}^{wt} \\
&= \sum_w \sum_k \sum_p \left( C_k^{wt} - C_p^{wt} \right) f_p^{wt} \rho_{pk}^{wt}.
\end{aligned} \quad (12)$$

Recall Property 1, we have $\sum_w \sum_k \sum_p \left( C_k^{wt} - C_p^{wt} \right) f_p^{wt} \rho_{pk}^{wt} \leq 0$. Hence, CNPSD is a RBAP. Replacing $\dot{f}_k^{wt}$ by $f_k^{w(t+1)} - f_k^{wt}$ above, we obtain the proof for DNPSD. □

**Corollary 2.** The stationary link flow pattern of NPSD is equivalent to Wardrop user equilibrium.

**Proof.** Zhang et al. (2001) had proved the stationary link flow pattern of RBAP is equivalent to Wardrop user equilibrium. According to Theorem 2, Corollary 2 holds. □

As is known to all, without extra behavioral conditions, Wardrop user equilibrium is not unique for path flow solution but only for link flow solution. However, according to Theorem 2 and Corollary 2, we can judge whether DNPSD has arrived at equilibrium by detecting the link flow states.

## 4 Lyapunov stability

Lyapunov method is an acknowledged methodology to analyze the stability of a dynamic system (Khalil 2002). The essence of this method is to identify a proper Lyapunov function. Peeta and Yang (2003) introduced the following Lyapunov function to analyze the stability of continuous-time PAP, namely

$$V(\mathbf{f}) = \int_0^{\Delta \mathbf{f}} \mathbf{c}(\mathbf{y}) d\mathbf{y}, \tag{13}$$

where $\Delta$ is the path link incidence matrix, the column vector $\mathbf{f}$ groups all the path flows and the row vector $\mathbf{c}(\cdot)$ groups the link travel times.

In this study, we still employ the above Lyapunov function to analyze the stability of the present CNPSD model.

**Theorem 3.** All the solutions of CNPSD are bounded and converge to Wardrop user equilibrium.

**Proof.** Since $\mathbf{c}(\cdot)$ is non-negative and $V(\mathbf{f}) \to \infty$ if $\mathbf{f} \to \infty$, then $V(\mathbf{f})$ is non-negative and radially unbounded. According to Eq. (12), take the derivative of $V(\mathbf{f})$ with respect to $t$, to give $\dot{V}(\mathbf{f}) = V'(\mathbf{f})\Delta \dot{\mathbf{f}} = \mathbf{c}\Delta \dot{\mathbf{f}} = \mathbf{C}\dot{\mathbf{f}} = \sum_w \sum_k C_k^{wt} \dot{f}_k^{wt} \leq 0$, where the row vector $\mathbf{C}$ groups all the route costs. Let $Z = \{\mathbf{f} | \dot{V}(\mathbf{f}) = 0\}$ and $E$ be the largest invariant set contained in the set $Z$, according to LaSalle's theorem (Khalil 2002), all the solutions of CNPSD model will be bounded and converge to the set $E$. Then, to prove Theorem 3, we only need to prove that the set $E$ is equivalent to Wardrop user equilibrium.

Since $\dot{V}(\mathbf{f}) = \sum_w \sum_k C_k^{wt} \dot{f}_k^{wt} = \sum_w \sum_k \sum_p \left(C_k^{wt} - C_p^{wt}\right) f_p^{wt} \rho_{pk}^{wt} \leq 0$, we have $\dot{V}(\mathbf{f}) = 0$ if and only if

$$\left(C_k^{wt} - C_p^{wt}\right) f_p^{wt} \rho_{pk}^{wt} = 0 \ \forall k, p, w. \tag{14}$$

According to Eq. (14), we have $C_p^{wt} \leq C_k^{wt}$ if $f_p^{wt} > 0$. Hence, given $f_p^{wt} > 0$ and $f_k^{wt} = 0$, we have $C_p^{wt} \leq C_k^{wt}$, which means that those used paths have smaller costs than those unused; again, given $f_p^{wt} > 0$ and $f_k^{wt} > 0$, we have $C_p^{wt} \geq C_k^{wt}$ and $C_p^{wt} \leq C_k^{wt}$, forcing $C_p^{wt} = C_k^{wt}$, which means that the used paths have the same cost. In summary, Eq. (14) implies that the used paths have the same cost, and this cost is smaller than the costs of the unused paths. Thus, $Z$ is exactly equivalent to Wardrop user equilibrium. From Property 3, the invariant set of CNPSD model is the feasible set $\Omega$ itself (defined in Eq. (10)). Thus, the largest invariant set $E = \Omega \cap Z = Z$. Then Theorem 3 is proved. □

Theorem 3 promises the convergence of CNPSD. However, we cannot conclude DNPSD is convergent from it. For DNPSD, we will examine the convergence through the numerical example conducted in Section 5.

## 5 Numerical example

In this section, using the day-to-day NPSD, we perform some numerical sensitivity analyses on an example network to interpret the impacts of reaction sensitivity on the traffic evolution and convergent characteristic of DNPSD. To this end, various levels of one-day capacity reductions are imposed on the example traffic network.

5.1 Example network and scenarios

The example network (see Fig. 2) has 12 nodes, 17 links, 2 OD-pairs and 8 paths, where OD-pair (1, 11) is connected by Path 1 ($1 \to 9 \to 14$), Path 2 ($1 \to 5 \to 10$), Path 3 ($2 \to 6 \to 10$) and Path 4 ($2 \to 11 \to 15$); OD-pair (2, 12) is connected by Path 5 ($3 \to 11 \to 16$), Path 6 ($3 \to 7 \to 12$), Path 7 ($4 \to 8 \to 12$) and Path 8 ($4 \to 13 \to 17$).

The travel demands for two OD-pairs are 90 (pcu/min). In Fig. 2, the bracketed three numbers in each link are the label, free-flow time (min) and normal capacity (pcu/min) of that link, respectively. The travel time of each link is computed by the BPR function

$$c_a^t = c_{a0}\left[1 + 0.15\left(\frac{v_a^t}{O_a}\right)^4\right] \quad \forall a, t, \tag{15}$$

where $c_{a0}$, and $O_a$ are the free-flow time and initial capacity of link $a$, respectively. For brevity, we will omit the units in the remaining context.

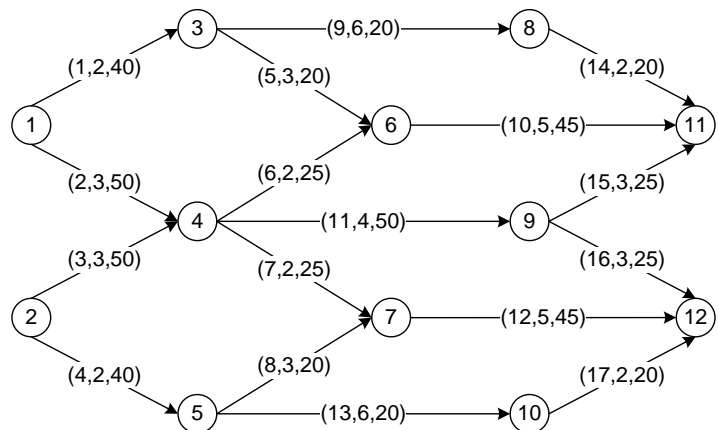

**Fig. 2** The example network

It can be concluded from Fig. 2 that the flows on Routes 1~4 respectively equal to that of Routes 8~5 when Fig. 2 is symmetric. This property will be often quoted below. Given that $s = 1.0E-10$, the initial traffic flow pattern of the network is at UE state with $\mathbf{f}^0 = (20, 20, 25, 25, 25, 25, 20, 20)^T$, and all travelers have the same reaction sensitivity $\theta$; let vector $\mathbf{Cap} = (Cap9, Cap11)$, where $Cap9$ and $Cap11$ are the reduction percentages of capacity on link 9 and link 11, respectively. In this numerical example, we design two scenarios which are parameterized as follows.

Asymmetric capacity reduction (ACR):
$\mathbf{Cap} = ([0.1:0.1:0.9], 0)$ when $\theta = [0.01:0.01:0.3]$.

Symmetric capacity reduction (SCR):
$\mathbf{Cap} = (0, [0.1:0.1:0.9])$ when $\theta = [0.01:0.01:0.3]$.

Note that the capacity reduces just at day 0 and reverts to normal after day 0. Also, only one link suffers from the capacity reduction at each time. Obviously, Fig. 2 remains symmetric under SCR, but becomes asymmetric under ACR. Here the convergence criterion is set to be $\left\| \mathbf{f}^{n+1} - \mathbf{f}^{n} \right\| \leq 1.0E - 005$.

5.2 Numerical results

Fig. 3 displays the traffic evolution processes under medium SCRs with $\mathbf{Cap} = (0, 50\%)$. For simplicity and without loss of generality, we just draw the results when $\theta = [0.05:0.05:0.3]$. Fig. 3 indicates the traffic evolutions on the symmetric routes are identical under all SCRs. It also indicates the traffic flows on the directly interfered routes (i.e., Route 4 and 5) face heavy flow losses, causing the sudden flow drops in the curves on the next day. In addition, when $\theta = [0.05:0.05:0.2]$, the route-swaps are convergent, whereas the smoothness of route-swaps becomes worse as $\theta$ rises. When $\theta = [0.25:0.05:0.3]$, the route-swaps fail to converge but finally reach the 2-day-cycled oscillations or quasi-oscillations, and the oscillation amplitudes increase as $\theta$ rises.

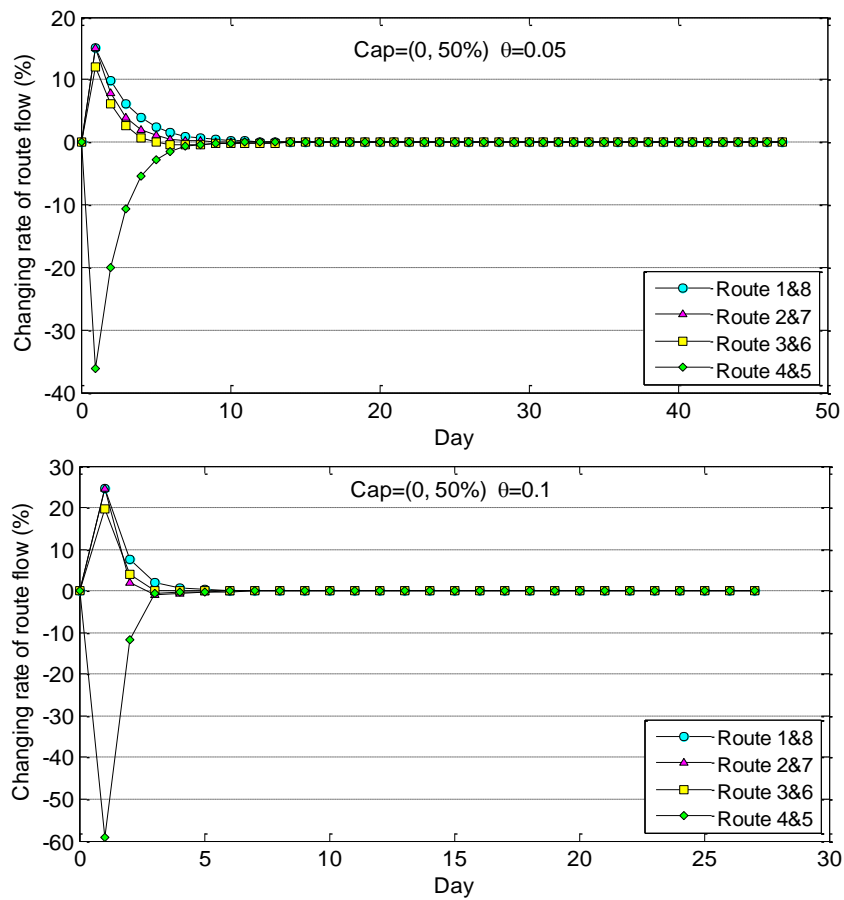

Published Online in *Networks and Spatial Economics* in Jan., 2015 (DOI: 10.1007/s11067-014-9281-3)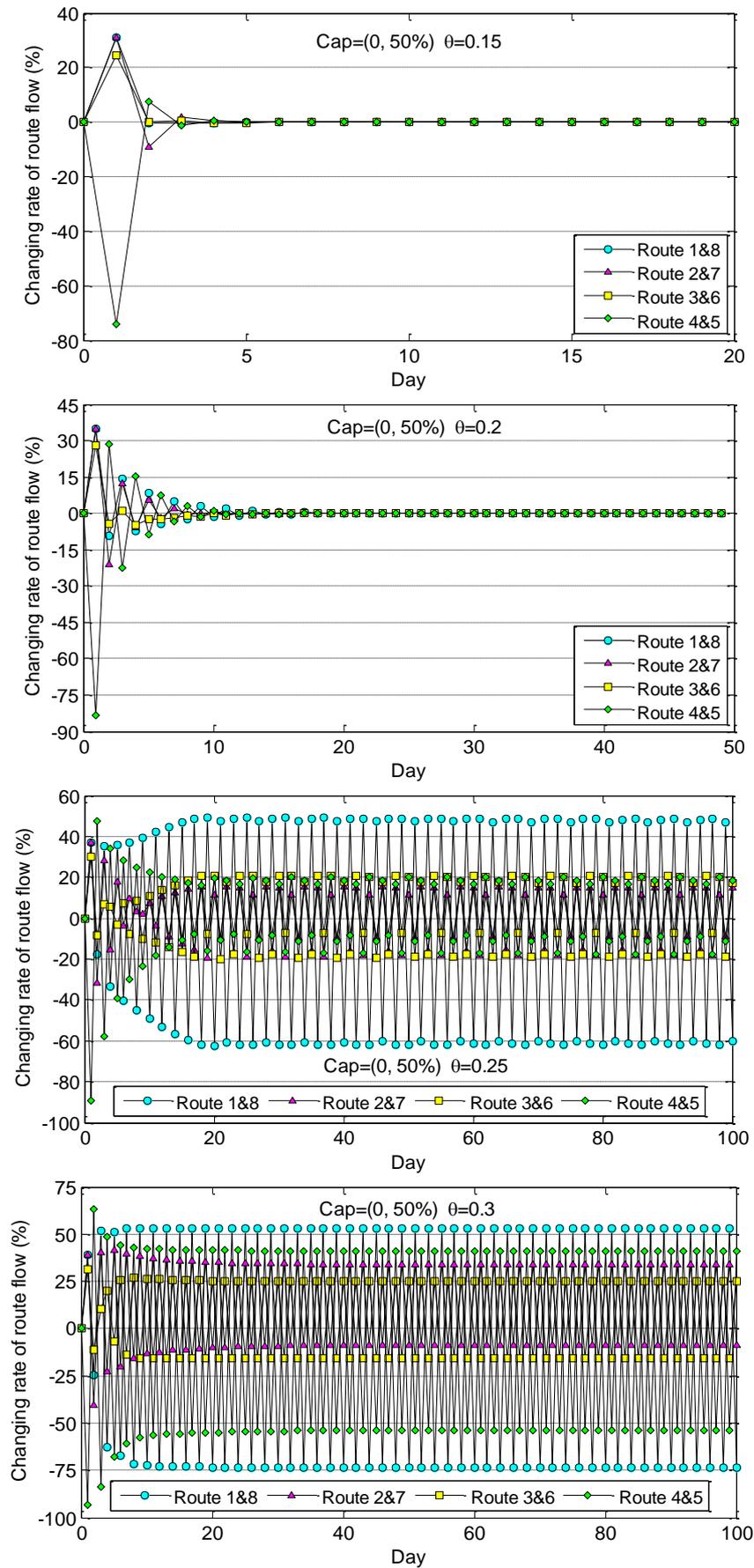

**Fig. 3** Traffic evolutions under different one-day medium SCRs

**13 / 19**

The detail numerical results show that if the route-swap is not able to converge, it will finally fall into 2-day-cycled oscillations or quasi-oscillations. To measure the final oscillation degrees, below we introduce the average deviation (AD) index, i.e.,

$$AD = \frac{1}{2}\sum_{t=1999}^{2000}\left[\sum_w \sum_k \left(f_k^{wt} - f_k^{w0}\right)^2\right]^{1/2}, \qquad (16)$$

where the number 2 in the denominator is the cycle of periodical (or quasi-periodical) oscillation. Since the route-swaps finally reach periodical (or quasi-periodical) oscillations within 100 days, here we draw the two path flow samples at iteration of 1999 and 2000. Of course, any two adjacent path flow samples after 100 days are effective. Obviously, AD is strictly positive if route-swaps are not convergent and 0 otherwise. In addition, AD increases when the oscillation becomes more severe. Subsequently, we present the ADs under two numerical scenarios.

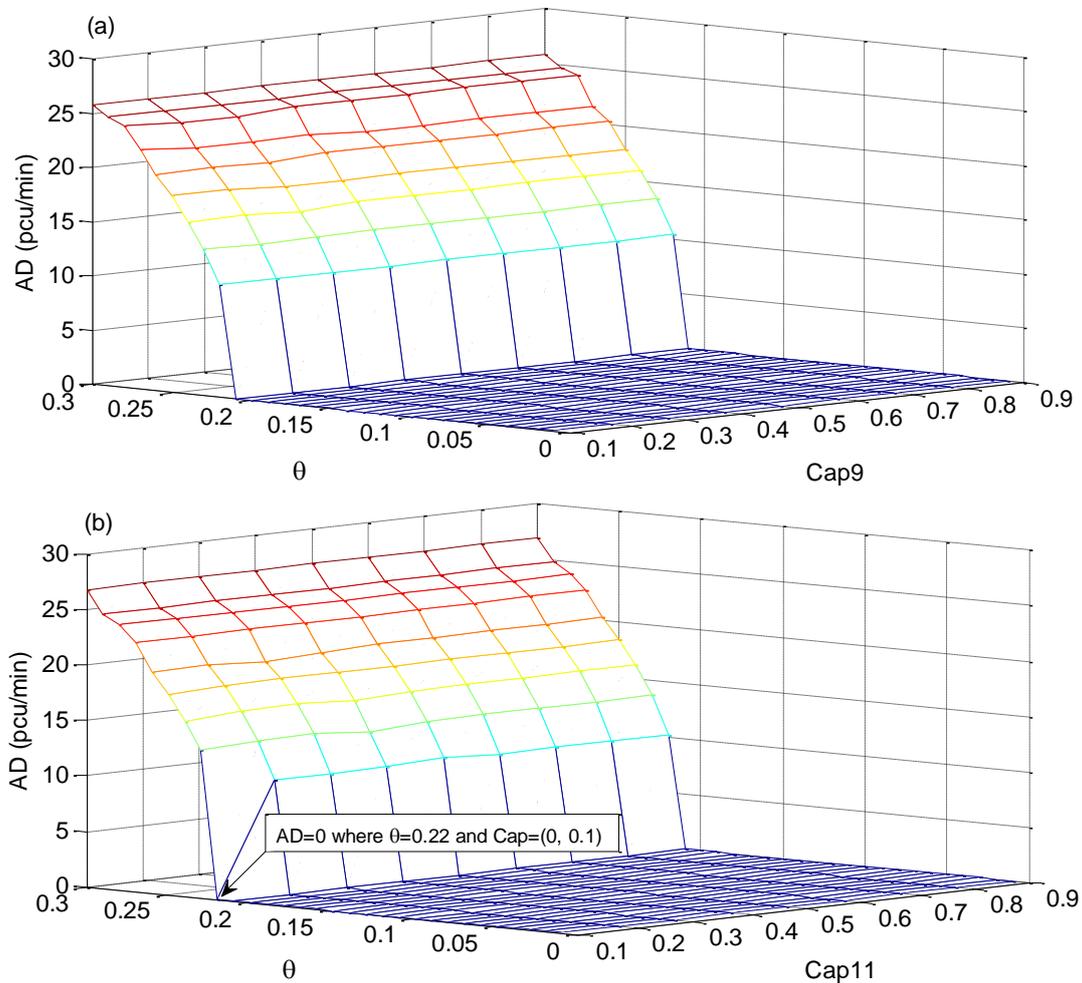

**Fig. 4** Final ADs under different $(\theta, \mathbf{Cap})$

Two sub-figures in Fig. 4 jointly shows that: 1) ADs are 0 when $\theta = [0.01:0.01:0.21]$ and then become larger than 0 and increase as $\theta$ increases in $\theta = [0.22:0.01:0.3]$; 2) for a specific $\theta$, ADs under different capacity reductions have few differences except the one when $\theta = 0$ and $\mathbf{Cap} = (0, 0.1)$. It demonstrates that the stability of DNPSD does not hold in general but influenced by the reaction sensitivity and capacity reduction; moreover, it depends much more on $\theta$ rather than the capacity reduction. Also, we can classify the final flow states into three phases, i.e., the stable phase with $\theta = [0.01 : 0.01 : 0.21]$, the meta-stable phase with $\theta = 0.22$ and the unstable phase with $\theta = [0.23 : 0.01 : 0.3]$. For the stable phase, the route-swaps can converge to the initial UE state for all capacity reductions; for the meta-stable phase, the route-swaps either converge or fall into 2-day-cycled quasi-oscillations in the end; for the unstable phase, the route-swaps finally fall into 2-day-cycled oscillations or quasi-oscillations.

Fig. 5 displays the impacts of capacity reduction $\mathbf{Cap}$ and reaction sensitivity $\theta$ on the convergent rate of route-swaps conducted by the day-to-day NPSD in the stable phase with $\theta = [0.01 : 0.01 : 0.21]$.

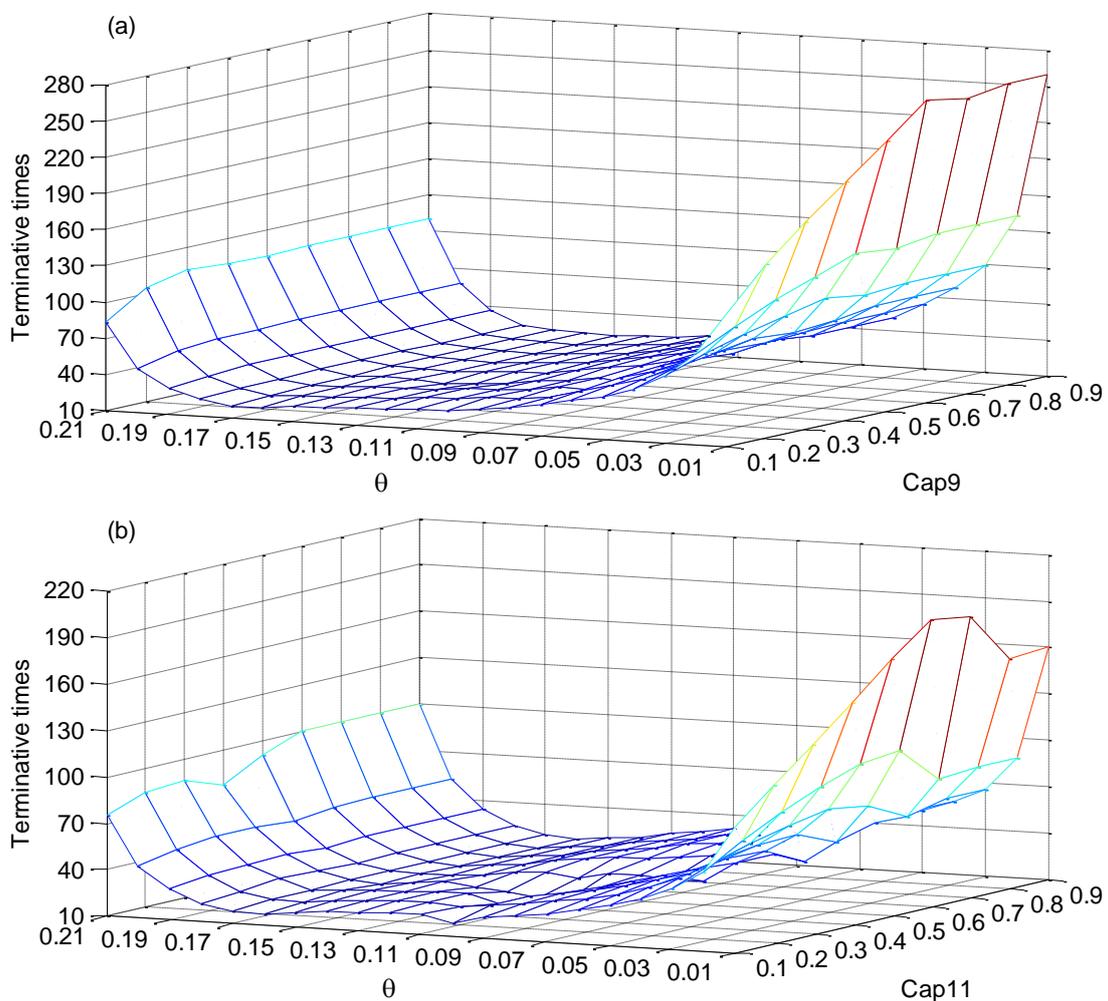

**Fig. 5** Convergence under different $(\theta, \mathbf{Cap})$

Fig. 5 suggests that, for a specific capacity reduction, the convergence rate increases at first and then decreases as $\theta$ increases from 0.01 to 0.3; for a given $\theta$, the convergence rate increases on the whole as the capacity reduction increases, but this trend weakens as $\theta$ increases. Overall, the reaction sensitivity $\theta$ has more significant effect on the convergent rate of route-swap than the capacity reduction **Cap**. The above findings can provide some useful insights for developing a NPSD-based swapping algorithm. Due to the value of $\theta$ cannot impact the feasibility of path flow solutions, we can maintain $\theta$ at a relative small (rather than a very small) level to realize a faster convergence. However, it should be noted that such a relative small $\theta$ may vary for different networks.

## 6 Conclusions and future researches

6.1 Research conclusions

In this paper, under the pairwise route-swapping behavioral framework, a nonlinear revision protocol is proposed and embedded into the traffic evolution equation of PAP, developing the present NPSD model. It is demonstrated that i) NPSD and PAP require the same amount of network information acquisition in route-swapping, ii) NPSD can prevent the over-swapping deficiency under a plausible behavior description, iii) NPSD can keep solution invariance, which means that the trial and error process for a feasible step-length can be omitted in a NPSD-based swapping algorithm and it helps raise the efficiency of algorithm, and iv) NPSD is a rational behavior swapping process and the CNPSD is globally convergent.

    Using the day-to-day NPSD, a numerical example is conducted to explore the effects of the reaction sensitivity on traffic evolution and the convergence characteristics of DNPSD. The following numerical results are found.

1) The final flow states can be divided into three phases, i.e., the stable phase, the meta-stable phase and the unstable phase. For the stable phase, route-swaps converge to UE under all designed capacity reductions; for the meta-stable phase, route-swaps either converge or fall into 2-day-cycled quasi-oscillations in the end; for the unstable phase, route-swaps finally fall into the 2-day-cycled oscillations or quasi-oscillations.
2) The reaction sensitivity $\theta$ has significant effect on the stability of DNPSD and the traffic evolution, while the capacity reduction takes a much smaller effect.
3) For DNPSD, there is a relative optimal interval associated with reaction sensitivity $\theta$ to realize a faster convergence.

6.2 Future researches

It should be noted that the numerical results stated in Section 6.1 are only applicable for the example network. Whether these results apply to the other ones needs to be investigated more exhaustively. Many further works are worthy of exploring based on NPSD.

1) From the behavioral aspect, empirical studies need to be performed to calculate $\theta$ and to understand the travelers' route-swapping behaviors better.
2) From the algorithmic aspect, i) the NPSD-based swapping algorithm deserves to be

developed (e.g., Zhang et al. 2014), compared with the other algorithms (e.g., Frank and Wolfe 1956; Huang and Lam 2002; Nie 2003; Patriksson 2004, etc.), and tested on the practical or large networks, and ii) a general stability theory for DNPSD needs to be developed.

**Acknowledgements** The authors extend their sincere thanks to the anonymous referees for their constructive comments. This study is jointly sponsored by the National Natural Science Foundation of China (71471014, 71401120) and the National Basic Research Program of China (2012CB725403-5).